\documentclass[11pt]{amsart}
\usepackage{amsmath,amsfonts,latexsym,graphicx,amssymb}
\usepackage{amsmath}

\theoremstyle{plain}
\newtheorem{theorem}{Theorem}[section]
\newtheorem{corollary}[theorem]{Corollary}
\newtheorem{lemma}[theorem]{Lemma}
\newtheorem{proposition}[theorem]{Proposition}

\def\sg {\mathrm{sg}}
\def\ur {\mathit{ur}}

\begin{document}   

\begin{abstract}
The Amitsur-Levitski theorem asserts that $M_n(F)$ satisfies a
polynomial identity of degree $2n$. (Here, $F$ is a field and
$M_n(F)$ is the algebra of $n \times n$ matrices over $F$). 
It is easy to give examples of subalgebras of $M_n(F)$ that do satisfy an identity of
lower degree and subalgebras of $M_n(F)$ that satisfy no polynomial identity of degree
$\le 2n-2$.
Our aim in this paper is to give a full classification of the subalgebras of
$n \times n$ matrices that satisfy no  nonzero polynomial of  degree less
than $2n$.
\end{abstract}

\title[On matrix subalgebras not satisfying an identity of degree $2n-2$.]
{On subalgebras of $n \times n$ matrices not satisfying identities of
  degree $2n-2$}
\author{Daniel Birmajer}
\address{Department of Mathematics \& Computer Science\\
    Nazareth College\\
4245 East Avenue\\Rochester, NY 14618-3790}

\email{birmajer@alumni.temple.edu}
  
\date{\today}   
\maketitle
\section{Introduction}\label{S:notation}
This paper is concerned with $n\times n$ matrix subalgebras that do not
satisfy a polynomial identity of degree $< 2n$.

To begin, let $F$ be a field,  $M_n(F)$ the algebra of $n \times n$ matrices over $F$,
and $F\left\{X\right\}=F\left\{X_1, X_2, \dots\right\}$ the free associative
algebra over $F$ in countably many variables.
A nonzero polynomial $f(X_1,\dots X_m) \in F\left\{X\right\}$ is a
\emph{polynomial identity} for an $F$-algebra $R$ (or, $R$ \emph{satisfies} $f$) if 
$f(r_1, \dots ,r_m)=0$ for all $r_1, \dots , r_m \in R$. 

Kaplansky (\cite{Ka48}) showed that if $R$ satisfies a polynomial of degree $d$, then it
satisfies a multilinear polynomial of degree $d$. The study of
identities for $R$ therefore reduces to the multilinear case.
The standard polynomial of degree $t$ is
\begin{equation*}
s_t(X_1, \dots, X_t) = \sum_{\sigma \in S_t} (\sg\sigma) 
X_{\sigma(1)}X_{\sigma(2)} \ldots X_{\sigma(t)},
\end{equation*} 
where $S_t$ is the symmetric group on $\{1, \dots, t\}$ and $(\sg\sigma)$ is the sign of 
the permutation $\sigma \in S_t$.
The standard polynomial $s_t$ is homogeneous of degree $t$, multilinear and alternating.
If $t$ is odd then $s_t(1,X_2,\dots , X_t)=s_{t-1}(X_2, \dots , X_t)$. Thus
$s_{2t}$ is an identity of $R$ if and only if $s_{2t+1}$ is an identity of $R$. 

The Amitsur-Levitski theorem asserts that $M_n(F)$ satisfies any standard
polynomial of degree $2n$ or higher. Moreover, if $M_n(F)$ satisfies a polynomial of degree
$2n$, then the polynomial is a scalar multiple of $s_{2n}$ (cf. \cite{AL50}).

The standard polynomial $s_{2n}$ is a minimal identity in the sense that
$M_n(F)$ satisfies no polynomial identity of degree less than $2n$. More
generally, if $A$ is a subalgebra 
of $M_n(F)$  isomorphic to a full block upper triangular matrix algebra,
\begin{equation*}
\begin{pmatrix}
 \fbox{$*$} &  &                 &    &   \\
   & \fbox{$*$}&              & \text{\huge *} &  \\
   &    & \ddots   &   &   \\
    & \text{\Large 0}   &               &   &  \fbox{$*$} \\
\end{pmatrix},
\end{equation*}
then $A$ satisfies no polynomial identity of  degree less than $2n$. To prove
this assertion, note that every full block upper triangular matrix algebra
contains the ``staircase sequence'' 
$e_{11}, e_{12}, e_{22}, e_{23}, \dots, e_{(n-1)(n-1)}, e_{(n-1)n}, e_{nn}$, and
\begin{equation}\label{E:staircase}
s_{2n-1}\left(e_{11}, e_{12}, e_{22}, e_{23}, \dots, e_{(n-1)(n-1)}, e_{(n-1)n}, e_{nn}\right)=e_{1n},
\end{equation}
where the $e_{ij}$ are the standard matrix units.

In $\S~2$ we provide the building blocks
for the main theorem of this paper and its proof. This proof and some of its
consequences are presented in  $\S~3$. For polynomial identities in ring theory
and the polynomial identities of $n\times n $ matrices,  \cite{Fo91} and
\cite{Ro80} are suggested general references.
\section{Building Blocks}\label{S:BB}
\begin{lemma}\label{L:Ai=Mi}
Let $A$ be a simple $F$-subalgebra of $M_n(F)$. Then either $A=M_n(F)$ or $A$
satisfies the identity $s_{2n-2}(A)=0$. 
\end{lemma}
\proof
By assumption, $A$ is a a finite dimensional central simple algebra over its center $k$.
Let $K$ denote the algebraic closure of $k$;  then   
$A \otimes_k K$ is a simple $K$-algebra in a natural way (cf. \cite{Ro80},\,\S 1.8),  with 
$\dim_K \left(A \otimes_k K\right) = \dim_k(A)$.
Also, $A \otimes_k  K \cong M_t(K)$ for some $t \le n$.
Suppose that  $A$ is a proper subalgebra of
$M_n(F)$. It follows that $t < n$.
Hence, by the Amitsur-Levitski theorem,  $A \otimes_k K$ satisfies $s_{2n-2}$,
and the result follows since $A$ is embedded as a
$k$ algebra in $A \otimes_k K$.
\endproof

Let $\ell , m$ be positive integers such that $\ell + m = n$ and set 
\begin{equation*}
E_{(\ell,m)} (F) = 
\begin{bmatrix}
 M_\ell(F) & M_{\ell \times m}(F)\\
 0         &  M_m(F)
\end{bmatrix},
\end{equation*}
an $F$-subalgebra of $M_n(F)$.

\begin{enumerate}
\item[(i)] Associated to  $E_{(\ell,m)}(F)$ are canonical F-algebra homomorphisms 
$$\pi_\ell\colon  E_{(\ell ,m)} (F) \to M_\ell (F) \text{\quad  and \quad} \pi_m\colon  E_{(\ell ,m)}(F) \rightarrow M_m(F).$$
Further identify $M_\ell(F)$ and $M_m(F)$ with 
\begin{equation*}
\begin{bmatrix}
 M_{\ell}(F)   & 0\\
 0             &  0
\end{bmatrix},
\begin{bmatrix}
 0         & 0\\
 0         &  M_{m}(F)
\end{bmatrix},
\end{equation*}
respectively.
\item[(ii)] Associated to a subalgebra $A$ of $E_{(\ell ,m)}(F)$ are
homomorphic image subalgebras $A_\ell$ and $A_m$ in $M_\ell (F)$ and 
$M_m (F)$ respectively.
\item[(iii)] Set 
\begin{equation*}
T_{(\ell,m)} (F) = 
\begin{bmatrix}
 0         & M_{\ell \times m}\\
 0         &  0
\end{bmatrix},
\end{equation*}
the Jacobson radical of $E_{(\ell ,m)}(F)$.
\end{enumerate}
\begin{lemma}\label{L:blocks}
Let $A$ be a subalgebra of $E_{(\ell ,m)} (F)$ such that $A_\ell$ satisfies 
$s_q$ for some $q \le 2\ell$ and  $A_m$ satisfies $s_r$ for some $r \le 2 m$. Then $A$ satisfies $s_{q+r}$.
\end{lemma}
\proof
Let $t=q+r$.
As an $F$-vector space, $E_{(\ell ,m)}(F) = M_\ell (F) \oplus T_{(\ell, m)} (F) \oplus  M_m (F)$. Thus each matrix $x$ in $A$ can be written as $x = a + b + c$ with $a \in A_\ell$, $b \in T_{(\ell,m)}$ and $c \in A_m$.  
Using linearity, we expand completely $s_t(x_1, \dots , x_t)$ and further use the following rules to simplify some of the terms:
\begin{enumerate}
\item $T_{(\ell,m)} (F)$ is a nilpotent ideal of $E_{(\ell ,m)}(F)$, with  $T^2_{(\ell,m)} (F) = 0$, 
and so each term in the expansion containing more than one entry in $T_{(\ell,m)} (F)$ equals 0. 
\item $M_\ell (F) M_m (F) = M_m (F) M_\ell (F)=0$.
\item $M_m (F) T_{(\ell,m)} (F) = T_{(\ell,m)} (F) M_\ell (F)=0$.
\end{enumerate}
We obtain
\begin{equation}\label{E:ac}
s_t(x_1, \dots ,x_n) =
\sum_{i=0}^{t+1}\sum_{\sigma \in S_t} (\sg\sigma) a_{\sigma(1)}\dots  a_{\sigma(i-1)}b_{\sigma(i)}c_{\sigma(i+1)}\dots  c_{\sigma(t)}.
\end{equation}
Fixing $i > q$, and given $\tau, \sigma \in S_t$, we say that $\tau$ is \emph{$i$-equivalent} to $\sigma$, if  $\tau$ restricted to the final interval $[i, t]$  equals the restriction of $\sigma$ to the same domain.
In symbols,
\begin{equation*}
\tau \sim_{i} \sigma \iff  \tau\vert_{[\,i,\,t]} = \sigma\vert_{[\,i,\,t]}.
\end{equation*}
For each $i>q$, the relation $\sim_{i}$ yields a partition of $S_t$ into disjoint subsets $P_{i}^k, \; k=1, \dots ,\frac{t!}{(i-1)!}$.
Then, we have
\begin{align*}
\sum_{\sigma \in S_t}& (\sg\sigma) a_{\sigma(1)}\dots  a_{\sigma(i-1)}b_{\sigma(i)}c_{\sigma(i+1)}\dots  c_{\sigma(t)} =\\
&= \sum_k \sum_{\sigma \in P_{i}^k} (\sg\sigma) a_{\sigma(1)}\dots  a_{\sigma(i-1)}b_{\sigma(i)}c_{\sigma(i+1)}\dots  c_{\sigma(t)}\\
&= \sum_k  (\sg \sigma_k) s_{i-1}(a_{\sigma_k(1)}, \dots , a_{\sigma_k(i-1)})b_{\sigma_k(i)}c_{\sigma_k(i+1)}\dots  c_{\sigma_k(t)},
\end{align*}
where  $\sigma_k$ is a representative of the class $P_{i}^k$. The last equality follows from the fact that for any $\sigma \in P_{i}^k$, $\sigma = \tau \circ \sigma_k$ for some $\tau  \in S_{i-1} \subseteq S_t$, and $(\sg \sigma) =  (\sg \tau)(\sg \sigma_k)$.
By assumption, $A_\ell$ satisfies $s_q$, and since $i-1 \ge q$ we obtain 
\begin{equation*}
\sum_{\sigma \in S_t}(\sg\sigma) a_{\sigma(1)}\dots  a_{\sigma(i-1)}b_{\sigma(i)}c_{\sigma(i+1)}\dots  c_{\sigma(t)} =0.
\end{equation*}
This shows that
\begin{equation}\label{E:a}
\sum_{i=q+1}^{t+1}\sum_{\sigma \in S_t} \sg(\sigma)\, a_{\sigma(1)}\dots  a_{\sigma(i-1)}b_{\sigma(i)}c_{\sigma(i+1)}\dots  c_{\sigma(t)}=0.
\end{equation}
For $i \le q$ we have that $t-i \ge r$. Applying a similar argument to the above, and using the fact that $A_m$ satisfies $s_r$, we see that also
\begin{equation}\label{E:c}
\sum_{i=0}^q\sum_{\sigma \in S_t} (\sg\sigma) a_{\sigma(1)}\dots  a_{\sigma(i-1)}b_{\sigma(i)}c_{\sigma(i+1)}\dots  c_{\sigma(t)}=0.
\end{equation}
Together, Equations \eqref{E:a} and \eqref{E:c} ensure that $s_t(x_1, \dots
,x_n) =0$, given Equation~\eqref{E:ac}.
\endproof
\subsection{} We now consider the case when $A$ contains a ``repetition''.
We will need some more notation.

\indent\textrm{(i)} Let $M_1, \dots M_t$ be matrices in $A$, 
\begin{equation*}
M_k =  \begin{bmatrix}
 a_k         & b_k    & c_k\\
 0             & e_k    & d_k\\
 0             & 0        & a_k
\end{bmatrix},\; a_k \in M_\ell (F), e_k \in M_m (F),  b_k \in M_{\ell \times m}(F), d_k \in M_{m \times \ell}(F).
\end{equation*}
Given $1 \le i < j \le t$ and $\sigma \in S_t$, set
\begin{equation*}
m_t^\sigma[i,j] = (\sg\sigma)\, a_{\sigma(1)}\dots  a_{\sigma(i-1)}b_{\sigma(i)}e_{\sigma(i+1)}\dots  e_{\sigma(j-1)}d_{\sigma(j)}a_{\sigma(j+1)}\dots  a_{\sigma(t)},
\end{equation*}
and denote by  $W$  the set of all matrix products
\begin{equation*}
\{ m_t^\sigma[i,j]: \sigma \in S_t \text{\,and\,} 1 \le i < j \le t\}. 
\end{equation*}
\indent\textrm{(ii)} The projection $\ur$ returns the $\ell \times \ell$ upper right block of a matrix in $A$:
\begin{equation*}
\ur\begin{bmatrix}
 a         & b    & c\\
 0         & e    & d\\
 0         & 0    & a
\end{bmatrix} = c
\end{equation*}
\indent\textrm{(iii)}
Given $n \times n$ matrices $M_1, \dots , M_t,$ we say that a matrix product  $M_1 \cdots M_t$
\emph{formally contains} the factor $A_1 \cdots A_s$ if $A_1=M_{\ell},
A_2=M_{\ell+1}, \dots ,A_s=M_{\ell+s-1}$, for some $1 \le \ell \le t$. This
notation is to distinguish to the case when $C A_1\cdots A_s D = M_1 \cdots
M_t$ as  $n \times n$ matrices, for some matrices $C$ and $D$.   
Further, if $\ell=1$, we say that $M_1 \cdots M_t$  formally contains $A_1
\cdots A_s$ as \emph{left factor}.

This is a good place to record a Lemma extracted from  \cite{AL50}, 
which will be used later.
\begin{lemma}\label{L:AL1}
\emph{\rm [\bf{AL50}, Lemma 1, 450-451]} 
If for an odd positive integer $r$ we put $Y=X_{i+1}\cdots X_{i+r},$ and if
$s'$ denotes the sum of all terms of $s_m(X)$ containing the
common factor $Y$, then 
\begin{equation*}
s'=s_{m-r+1}(X_1, \dots , X_i, Y, X_{i+r+1}, \dots , X_m).
\end{equation*}
\end{lemma} 
\begin{lemma}\label{L:ur}
Set $t = 2 (\ell + m),$ and let $M_1, \dots, M_t$ be matrices in $A$ such that
for all $1 \le k \le t$, 
\begin{equation*}
M_k =  \begin{bmatrix}
 a_k         & b_k    & 0\\
 0             & e_k    & d_k\\
 0             & 0        & a_k
\end{bmatrix},\text{\; for \;} a_k \in M_\ell (F), e_k \in M_m (F),  b_k \in M_{\ell \times m}(F), d_k \in M_{m \times \ell}(F).
\end{equation*}
Then $\ur\left[s_t(M_1, \dots ,M_t)\right] = 0$.
\end{lemma}
\proof
First we observe that
\begin{equation*}
\ur[M_1 \cdots M_t] = \sum_{1 \le i < j \le t} a_{\sigma(1)}\dots  a_{\sigma(i-1)}b_{\sigma(i)}e_{\sigma(i+1)}\dots  e_{\sigma(j-1)}d_{\sigma(j)}a_{\sigma(j+1)}\dots  a_{\sigma(t)},
\end{equation*}
which implies that 
\begin{equation}\label{E:ur}
\ur\left[s_t(M_1, \dots ,M_t)\right] = \sum_{\sigma \in S_t}\; \sum_{1 \le i < j \le t} m_t^\sigma[i,j].
\end{equation}
To prove that $\ur\left[s_t(M_1, \dots ,M_t)\right] = 0$, we split the right hand side into two summands:
\begin{align}
\ur &\left[s_t(M_1, \dots ,M_t)\right] =\notag\\
&\sum_{\sigma \in S_t} \sum_{\begin{subarray}{1}1 \le i < j \le t\\
\\j-i-1 \ge 2m\end{subarray}}  m_t^\sigma[i,j] + 
\sum_{\sigma \in S_t} \sum_{\begin{subarray}{1}1 \le i < j \le t\\
\\ j-i \le 2m \end{subarray}}  m_t^\sigma[i,j]\label{E:7}
\end{align}
Our goal is to show that each summand in \eqref{E:7} is zero. To handle the
first summand we introduce the following new equivalence relation on $S_t$. 
Given fixed $1 \le i < j \le t$, such that $j-i-1 \ge 2m$, and given $\tau,
\sigma \in S_t$, say that $\tau$ is 
\emph{$[i,j]$-equivalent} to $\sigma$ if  $\tau$ restricted to the initial
and final intervals 
$[1,i]$ and $[j,t]$ equals the restriction of $\sigma$ to the same domain.
In symbols,
\begin{equation*}
\tau \sim_{[i,j]} \sigma \iff \tau\vert_{[1,i]} = \sigma\vert_{ [1,i]} \text{\; and\;} \tau\vert_{ [j,t]} = \sigma\vert_{ [j,t]}
\end{equation*}
For each pair $i,j$, such that $1 \le i < j \le t$ and $j-i-1 \ge 2m$, the
relation $\sim_{[i,j]}$ yields a partition of 
$S_t$ into disjoint subsets $P_{[i,j]}^k, \; k=1, \dots ,
\frac{t!}{(j-i-1)!}$.
Then, we have
\begin{align*}
&\sum_{\sigma \in S_t}\; \sum_{\begin{subarray}{1}1 \le i < j \le t\\
\\j-i-1 \ge 2m\end{subarray}}  m_t^\sigma[i,j] =
\sum_{\begin{subarray}{1}1 \le i < j \le t\\ \\j-i-1 \ge 2m\end{subarray}} 
\sum_k \sum_{\sigma \in P_{[i,j]}^k}  m_t^\sigma[i,j]=\\\\
&\sum_{\begin{subarray}{1}1 \le i < j \le t\\ \\j-i-1 \ge 2m\end{subarray}} 
\sum_k \sum_{\sigma \in P_{[i,j]}^k} (\sg\sigma)\, a_{\sigma(1)}\cdots
a_{\sigma(i-1)}b_{\sigma(i)}e_{\sigma(i+1)}\cdots
e_{\sigma(j-1)}d_{\sigma(j)}a_{\sigma(j+1)}\cdots  a_{\sigma(t)}\\
&\qquad =\\
&\sum_{\begin{subarray}{1}1 \le i < j \le t\\ \\j-i-1 \ge 2m\end{subarray}} 
\sum_k (\sg \sigma_k) \, a_{\sigma_k(1)}\cdots
 a_{\sigma_k(i-1)}b_{\sigma_k(i)}\,s\, d_{\sigma_k(j)}a_{\sigma_k(j+1)}\cdots  a_{\sigma_k(t)},
\end{align*}
where  $s=s_{i-j+1}(e_{\sigma_k(i+1)},\dots , 
 e_{\sigma_k(j-1)})$ and $\sigma_k$ is a representative of the class $P_{[i,j]}^k$.
Since $j-i-1 \ge 2m$, 
\begin{equation*}
s_{i-j+1}(e_{\sigma_k(i+1)},\dots , e_{\sigma_k(j-1)})=0
\text{\quad for all\;} k,
\end{equation*}
hence
\begin{equation*}
\sum_{\sigma \in S_t} \;\sum_{\begin{subarray}{1}1 \le i < j \le t\\
\\j-i-1 \ge 2m\end{subarray}}  m_t^\sigma[i,j] =0.
\end{equation*}
This takes care of the first term in  \eqref{E:7}.
We now turn to the second summand.
For a given $q$, with $2 \le q \le t$, denote by $R_q$ the set of all
 $q$-tuples $r = (r_1, \dots , r_q)$ of different
elements from $\{1,\dots ,t\}$ and by $T_{(r_1, \dots , r_q)}$ the set of
matrix products $w$ formally containing the common factor 
$b_{r_1} e_{r_2}\cdots e_{r_{q-1}}d_{r_q}$.
Considering all possible $q$ and $q$-tuples, the sets $T_{(r_1, \dots , r_q)}$
form a partition of $W$. 
We are interested in the case when $q \le 2m+1$. Observe that
\begin{equation*}
\sum_{\sigma \in S_t}\; \sum_{\begin{subarray}{1}1 \le i < j \le t\\ \\j-i \le 2m\end{subarray}}   m_t^\sigma[i,j]
=\sum _{q=2}^{2m+1}\; \sum_{r \in R_q} \;\sum_{w \in T_{(r_1, \dots , r_q)}} w.
\end{equation*}
Fix $q$ odd, a $q$-tuple $(r_1, \dots , r_q)$, and the corresponding set of
matrix products $T_{(r_1, \dots , r_q)}$.
Then, $\sum_{w \in T_{(r_1, \dots , r_q)}} w$ is the sum of all matrix products
formally containing the common factor 
$y = b_{r_1} e_{r_2}\cdots e_{r_{q-1}}d_{r_q}$.
Each matrix product $w  \in T_{(r_1, \dots , r_q)}$ corresponds uniquely to a
permutation $\sigma \in S_t$ and a  pair 
$(i,j)$, such that  the $q$-tuple $(r_1, \dots ,
r_q)$ is the image under $\sigma$ of 
$(i, \dots ,j)$. Explicitely, the correspondence is  $w= m_t^\sigma[i,j]$.
We can now apply  Lemma~ \ref{L:AL1} and the alternating property of the
standard polynomials. If $\sigma_0 \in S_t$ is a fixed permutation such that 
$$\sigma_0\colon i \rightarrow r_i, \, 1 \le i \le q,$$  we have
\begin{equation*} 
\sum_{w \in T_{(r_1, \dots , r_q)}} w
=(\sg\, \sigma_0)\, s_{t-q+1}\left(y, a_{\sigma_0(q+1)}, \dots , a_{\sigma_0(t)}\right)
\end{equation*}
where  $y = b_{r_1} e_{r_2}\cdots e_{r_{q-1}}d_{r_q}$. 
Since $t-q+1 \ge 2\ell$, and since all the arguments of  $s_{t-q+1}$ in the last
equation are $\ell \times \ell$ matrices, it 
follows that
\begin{equation}\label{E:qodd}
\sum_{w \in T_{(r_1, \dots , r_q)}} w = 0, \; \text{when}\; q \; 
\text{is odd and $(r_1, \dots , r_q)$ is a fixed $q$-tuple.}
\end{equation}
Therefore
\begin{equation*}
\sum_{\begin{subarray}{1}q=2\\\text{$q$ odd}\end{subarray}}^{2m+1}
\; \sum_{r \in R_q} \;\sum_{w \in T_{(r_1, \dots , r_q)}} w = 0.
\end{equation*}
Suppose now that  $q$ is even, so $q \le 2m$, and fix an arbitrary $q$-tuple
$r= (r_1, \dots , r_q)$. 
We will split further the  sets $T_r$. First consider all $w \in T_r$ formally
containing in common the left factor 
$y = b_{r_1} e_{r_2}\cdots e_{r_{q-1}}d_{r_q}$, and call this subset  $L_r$. 
Then, for each $r_0  \not\in \{r_1, \dots ,r_q\}$ consider the $(q+1)$-tuple $(r_0,
r)$ and  the subset 
$G_{(r_0, r)}$ of $w \in T_r$ formally containing in common the factor $y=a_{r_0}b_{r_1} e_{r_2}\cdots e_{r_{q-1}}d_{r_q}$. 
The sum of all matrix products in the set $T_r$ can be split as
\begin{equation*}
\sum_{w \in T_r} w = \sum_{w \in L_r}w
+ \sum_{r_0: r_0 \ne r_1, \dots ,r_q}\;\; \sum_{w \in G_{(r_0,r)}}w.
\end{equation*}
For the terms in $L_r$ we have
\begin{equation}
\sum_{w \in L_{(r_1, \dots , r_q)}} w
=(\sg\, \sigma_0) \, y \,  s_{t-q}\left(a_{\sigma_0(q+1)}, \dots , a_{\sigma_0(t)}\right),
\end{equation}
where  $y = b_{r_1} e_{r_2}\cdots e_{r_{q-1}}d_{r_q}$,  and where $\sigma_0
\in S_t$ is a fixed permutation such that 
$$\sigma_0\colon i \rightarrow r_i, \, 1 \le i \le q.$$
Since $t-q \ge 2\ell$, we obtain 
\begin{equation}
\sum_{w \in L_r}  w = 0.
\end{equation}

Finally, for a suitable fixed $r_0$, the sequence $(r_0,r)$  has odd length, so we can argue as in~\eqref{E:qodd} to obtain   
\begin{equation*}
\sum_{w \in G_{(r_0,r)}} w
=(\sg\, \sigma_0) \, s_{t-q+1}\left(y, a_{\sigma_0(q+2)}, \dots , a_{\sigma_0(t)}\right) = 0,
\end{equation*}
where $y = a_{r_0}b_{r_1} e_{r_2}\cdots e_{r_{q-1}}d_{r_q}$, 
and where $\sigma_0 \in S_t$ is a fixed permutation such that 
\begin{equation*}
\sigma_0=
\begin{cases}
1\rightarrow r_0,\\
i \rightarrow r_{i-1},\text{\quad for \quad} 2 \le i \le q+1.
\end{cases}
\end{equation*}
This finishes the proof of Lemma~\ref{L:ur}.
\endproof
\begin{proposition}\label{P:abced}
Let 
\begin{equation*}
A =  \left\{\begin{bmatrix}
 a         & b    & c\\
 0         & e    & d\\
 0         & 0    & a
\end{bmatrix}:
a, c \in M_\ell (F), e \in M_m (F),  b \in M_{\ell \times m}(F), d \in M_{m \times \ell}(F)  \right\}.
\end{equation*}
Then, $A$ satisfies $s_{2(\ell + m)}$.
\end{proposition}

\proof For any $t$ and matrices $M_k \in A, \, k = 1 \dots t$, set 
\begin{equation*}
M_k = \begin{bmatrix}
 a_k         & b_k    & c_k\\
 0             & e_k    & d_k\\
 0             & 0        & a_k
\end{bmatrix},\;
a_k \in M_\ell (F), e_k \in M_m (F),  b_k \in M_{\ell \times m}(F), d_k \in M_{m \times \ell}(F).
\end{equation*}
By direct calculations, we obtain 
\begin{align*}
\ur&[s_t(M_1, \dots , M_t)] =\\
&=\sum_{i=1}^t s_t(a_1, \dots, a_{i-1}, c_i,a_{i+1}, \dots , a_t)
+ \sum_{\sigma \in S_t}\; \sum_{1 \le i < j \le t}
 m_t^\sigma[i,j].
\end{align*} 
Now set $t=2(\ell+m)$. It follows from \eqref{E:ur} that
\begin{equation*}
 \sum_{\sigma \in S_t} \;\sum_{1 \le i < j \le t}
 m_t^\sigma[i,j]=\ur\left[s_t(M_1^\prime, \dots ,M_t^\prime)\right]=0,
\end{equation*}
where $M_k^\prime$ is the matrix in $A$ obtained by replacing the upper right corner
$c_k$ of $M_k$ by $0 \in M_\ell(F)$.
Suitable applications of the  Amitsur-Levitski identity give us
\begin{equation*}
\ur\left[s_t(M_1, \dots , M_t)\right]=0,
\end{equation*}
\begin{equation*}
s_t\left(
\begin{bmatrix}
 a_1         & b_1\\
 0           & e_1    
\end{bmatrix}, \dots ,
\begin{bmatrix}
 a_t         & b_t\\
 0           & e_t    
\end{bmatrix}\right)=0,
\end{equation*}
and 
\begin{equation*}
s_t\left(
\begin{bmatrix}
 e_1         & d_1\\
 0           & a_1    
\end{bmatrix}, \dots ,
\begin{bmatrix}
 e_t         & d_t\\
 0           & a_t    
\end{bmatrix}\right)=0.
\end{equation*}
Combining the three equations, it follows that $s_t\left(M_1, \dots ,  M_t\right)=0.$
\endproof
\section{Main Theorem}\label{S:main}
In this section we prove that if a matrix subalgebra of $M_n(F)$ does not
satisfy the standard polynomial 
$s_{2n-2}$, then it is isomorphic as $F$-algebra to a full block upper
triangular matrix algebra.
\subsection{} We first introduce our notation and review some necessary background (cf. \cite{Le02}).

\indent\textrm{(i)}\;
Let $t$ be a positive integer, let $\ell_1,\ell_2  ,\cdots ,\ell_t$ be positive integers summing up to $n$, and set 
\begin{equation*}
E_{(\ell_1,\ell_2  ,\dots ,\ell_t)} (F) = 
\begin{bmatrix}
M_{\ell_1}(F) & M_{\ell_1 \times \ell_2}(F) & \cdots & M_{\ell_{1} \times \ell_{t-1}}(F) & M_{\ell_{1} \times \ell_t}(F) \\
 0         &  M_{\ell_2}(F) & \cdots & M_{\ell_{2} \times \ell_{t-1}}(F) & M_{\ell_{2} \times \ell_t}(F)\\
\vdots        & \vdots       &\vdots                          &\vdots            &\vdots                \\
0    & 0  & \cdots    &  M_{\ell_{t-1}}(F)  & M_{\ell_{t-1} \times \ell_t}(F) \\
0    & 0  & \cdots    & 0                          & M_{\ell_t}(F) \\
\end{bmatrix},
\end{equation*}
a full block upper triangular matrix subalgebra of $M_n(F)$.

\indent\textrm{(ii)}\; Recall that every $F$-algebra automorphism $\tau$ of $M_n(F)$ 
is inner (i.e., there exists an invertible $Q$ in $M_n(F)$ 
such that $\tau(a)=QaQ^{-1}$ for all $a \in M_n(F)).$
We will say that two $F$-subalgebras $A,\, A^{\prime}$ of $M_n(F)$ are
\emph{equivalent} provided there exists an automorphism 
$\tau$ of $M_n(F)$ such that $\tau(A) = A^{\prime}$.

\indent\textrm{(iii)}\; Associated to $E_{(\ell_1,\ell_2  ,\dots ,\ell_t)}
(F)$ are canonical F-algebra homomorphisms 
\begin{equation*}
\pi_{ij}:E_{(\ell_1,\ell_2  ,\dots ,\ell_t)} (F) \to  E_{(\ell_i, \ell_{i+1},
  \dots ,\ell_j)} (F), \text{\; for \;} 1 \le i \le j \le t.
\end{equation*}
When $i=j$ we write $\pi_i$ for $\pi_{ii}$.
For a subalgebra $\Lambda$ of $E_{(\ell_1,\ell_2  ,\dots ,\ell_t)} (F)$, we have the homomorphic images:  
\begin{equation*}
\Lambda_{ij} := \pi_{ij}(\Lambda),
\end{equation*}
embedded in $E_{(\ell_i ,\ell_{i+1} ,\dots ,\ell_j)}$.

\indent\textrm{(iv)}\;
We will say that a subalgebra $\Lambda$ of  $E_{(\ell_1,\ell_2,\dots ,\ell_t)} (F)$ is an
 \emph{$(\ell_1,\ell_2  ,\dots ,\ell_t)$-extension of simple blocks} if the restrictions 
$\pi_i:\Lambda \to M_{\ell_i}(F)$, for $1\le i \le t$, are all irreducible
representations (when $F$ is algebraically closed, of course, the
representation $\pi_i$ is irreducible if and only if
$\pi_i(\Lambda)=M_{\ell_i}$). Note that, every $F$-subalgebra $A$ of $M_n(F)$ is equivalent to an
$(\ell_1,\ell_2  ,\dots ,\ell_t)$-extension of 
simple blocks $\Lambda$ for some suitable $(\ell_1,\ell_2  ,\dots ,\ell_t)$.

\indent\textrm{(v)}\;
Further we will say that $\Lambda$  \emph{contains a repetition} when
$$\pi_i:\Lambda \to M_{\ell_i} \text{\quad  and \quad}   \pi_j:\Lambda \to
M_{\ell_j}$$
 are equivalent representations, for some $1 \le i < j \le t$ (and so $\ell_i = \ell_j$).
Also, $\Lambda$ is \emph{uniserial} when $\Lambda_{i (i+1)}$ is not semisimple, for all $1 \le i \le (t-1)$.

\begin{lemma}\label{L:repetition}
If an extension of simple blocks $\Lambda$ contains a repetition, then the
 standard identity $s_{2n-2}=0$ holds
for $\Lambda$. 
\end{lemma}
\proof
Assume $\pi_i:\Lambda \to M_{\ell_i}$ and   $\pi_j:\Lambda \to M_{\ell_j}$ are equivalent representations for some $1 \le i < j \le t$. Then we can choose an $F$-algebra automorphism $\tau$ of $M_n(F)$ such that $\pi_{ij}(\tau(\Lambda))$ is a subalgebra of  
\begin{equation*}
 \left\{\begin{bmatrix}
 a         & b    & c\\
 0         & e    & d\\
 0         & 0    & a
\end{bmatrix}:
a, c \in M_{\ell_i} (F), e \in M_{\ell_{i+1}} (F),  b \in M_{\ell_i \times \ell_{i+1}}(F), d \in M_{ \ell_{i+1} \times \ell_i}(F)  \right\}.
\end{equation*}
The result now follows from Lemma~\ref{L:blocks} and Proposition~\ref{P:abced}. 
\endproof 
\begin{lemma}\label{L:simple}
If an extension of simple blocks $\Lambda$ is not uniserial, then the standard identity $s_{2n-2}=0$ holds for $\Lambda$. 
\end{lemma}
\proof
Follows immediately from Lemma~\ref{L:blocks}.
\endproof
\begin{theorem}\label{T:main} 
Let $F$ be a field and let $A$ be an $F$-subalgebra of $M_n(F)$. If $A$ does not satisfy the standard polynomial $s_{2n-2}$, then $A$ is equivalent to a full block upper triangular matrix algebra. 
\end{theorem}
\proof
It suffices to show that the only 
$(\ell_1,\ell_2  ,\dots,\ell_t)$-extension of simple blocks 
$\Lambda$ for which  the standard polynomial $s_{2n-2}$ is not an identity is
the full block upper triangular 
matrix algebra $E_{(\ell_1,\ell_2  ,\dots ,\ell_t)} (F)$.   
By Lemma~\ref{L:Ai=Mi}, $\Lambda_i = M_{\ell_i}(F)$ for $1 \le i \le t$.
By Lemma~\ref{L:simple} and Lemma~\ref{L:repetition}, $\Lambda_{i(i+1)}(F)$ is not semisimple
and does not contain a repetition, for each $1 \le i \le t-1$.   
We  conclude that (cf. \cite{Le02}, Lemma 3.6)
\begin{equation*}
\Lambda_{i(i+1)}(F) = M_{\ell_i \times \ell_{i+1}}(F), \text{\quad for
  each\quad} 1 \le i \le t-1.   
\end{equation*}
Therefore, $\Lambda$ contains the
staircase unit matrices (c.f. \eqref{E:staircase}), and  every unit matrix
$e_{ij}, \text{\; for \;} j > i\;$ can be expressed
as a product of those. The Theorem now follows.
\endproof
\begin{corollary}
The standard polynomial $s_{2n-2}$ is an identity for any proper subalgebra of
$U_n(F)$, the algebra of upper triangular matrices over the field $F$.
\end{corollary}
\proof
Immediate from Theorem~\ref{T:main}.
\endproof

\noindent
\textbf{Remark}
The standard polynomial of degree $2n-2$ is not necessarily an identity for any proper
subalgebra of $U_n(C)$ when $C$ is a commutative ring: Let $I$ be a nonzero
ideal of $C$, and consider the $C$-subalgebra $B$ of $U_n(C)$ defined by the
property that the $(1,2)$-entry of matrices in $B$ lie in $I$. A staircase
argument shows that $s_{2n-2}(B) \ne 0$.
\newpage
\section*{Acknowledgments}
The results in this paper are part of the author's Ph.D. thesis at Temple University, and the
author thanks his advisor Edward Letzter for his help and guidance.

\enddocument